\documentstyle{article}
\newcommand{\bt}{\begin{tabular}}
\newcommand{\et}{\end{tabular}}
\newcommand{\ba}{\begin{array}}
\newcommand{\ea}{\end{array}}
\newcommand{\be}{\begin{equation}}
\newcommand{\ee}{\end{equation}}
\newcommand{\ben}{\begin{enumerate}}
\newcommand{\een}{\end{enumerate}}

\newcommand{\Ttimes}{\mbox{$\hspace{.5mm}\bigcirc\hspace{-3.2mm}
\raisebox{-.7mm}{$\top$}\hspace{1mm}$}}

\newtheorem{Thm}{Theorem}[section]
\newtheorem{Prop}[Thm]{Proposition}

\newtheorem{Rem}[Thm]{Remark}

\newcommand{\dowod}{\noindent{\bf Proof:} }
\newcommand{\qed}{\hfill $\Box $ }

\newcommand{\End}{\mbox{\rm End}\, }

\newcommand{\Ad}{\mbox{\rm Ad}\, }
\newcommand{\id}{\mbox{\rm id}\, }
\newcommand{\Ldwa}{\,\stackrel{2}{\bigwedge}}
\newcommand{\Odwa}{\,\stackrel{2}{\bigotimes}}

\newcommand{\w}{{\!}\wedge{\!}}

\newcommand{\pib}{\pi _{\Join}}
\def\sfR{\mbox{\bf R} }

\font \eul=eufm10 scaled \magstep2

\newcommand{\gotG}{\mbox{\eul g}}

\newcommand{\er}{\varepsilon }

\newcommand{\sir}{\sigma }
\def\Dr{\Delta }

\def\Drop{\Delta ^{\rm op}}

\begin{document}

\title{\bf Gauge transformations and quasitriangularity}
\author{{\bf S. Zakrzewski}  \\
\small{Department of Mathematical Methods in Physics,
University of Warsaw} \\ \small{Ho\.{z}a 74, 00-682 Warsaw, Poland} }
\date{}
\maketitle

\begin{abstract}
Natural conditions on a Poisson/quantum group $G$ 
to implement Poisson/quantum gauge transformations on the lattice
are investigated. In addition to our previous result that
transformations on one lattice link require $G$ to be coboundary,
it is shown that for a sequence of links one needs 
a quasitriangular $G$. 
\end{abstract}

\section{Introduction: one link results}

Is it possible to formulate lattice (group valued) gauge fields
with a Poisson (or quantum) group as the structure group, in such a
way that the gauge transformations are Poisson or quantum actions?
Which conditions should such a group satisfy?

In \cite{cob} we have answered the last question in case of
gauge transformations performed on one link only. We recall
these results.
\ben
\item Let $(G,\pi )$ be a Poisson Lie group. There exists a
Poisson structure $\sigma $ on $G$ such that the map
\be\label{1lin}
 G\times G\times G \ni (x,y,z)\mapsto xyz^{-1}\in G
\ee
is Poisson as a map from $(G,\pi)\times (G,\sigma )\times (G,\pi
)$ to $(G,\sigma )$ {\em if and only if} $(G,\pi )$ is coboundary.
In this case $\sigma $ is just the `plus' structure
\be
\sigma (g) = \pi _+ (g):= rg+gr,
\ee
where  $r\in \Ldwa \gotG$ such that $\pi (g)=rg-gr$.
\item Let $(H,m,\Dr )$ be a Hopf algebra. Let $\Drop$ denote the
comultiplication opposite to $\Dr $: $\Drop = P\circ\Dr$, where
$P$ is the permutation in the tensor product. 

There exists a (new)
coalgebra structure $\Sigma\colon H\to H\otimes H$ such that the map
\be
m(m\otimes m)\colon H\otimes H\otimes H\to H
\ee
is a morphism from $(H,\Dr )\otimes (H,\Sigma )\otimes (H,\Drop )$
to $(H,\Sigma )$ {\em if and only if} there exists $\sfR\in H\otimes
H$ such that
\be 
\Dr (a)\sfR = \sfR \Drop (a)\qquad \mbox{\rm and}\qquad
[(\Dr \otimes \id )\,\sfR  \,]\,\sfR _{12}=[(\id\otimes \Dr )\,\sfR
\, ] \sfR _{23}.
\ee
In this case, $\Sigma $ is given by the `plus structure'
\be
\Sigma (a) = \Dr _+ (a):= \Dr (a)\, \sfR .
\ee
\item For matrix quantum groups with basic relations given by 
\be
\hat{R}(u\Ttimes u)=(u\Ttimes u)\hat{R},
\ee
the basic `plus' commutation relations are given by
\be
\hat{R}(u\Ttimes u)=(u\Ttimes u)\hat{R}_{21}\qquad\qquad
(R_{21}\equiv P\hat{R}P). 
\ee
\een
\section{Two links --- the results}

In this paper we consider the case of gauge transformations on
two consecutive links of the lattice. In the sequel, we
concentrate mainly on the Poisson case. We assume already, that we
have the `plus' structure on links, and we investigate conditions
under which the map
\be\label{2lin}
G\times G\times G \ni (a,g,b)\mapsto (ag^{-1},gb)\in G\times G
\ee
is Poisson as a map from $(G,\pi _+)\times (G,\pi )\times (G,\pi
_+)$ to $(G,\pi _+)\times (G,\pi _+)$. 
Due to appearance of the same $g$ in both components of map
(\ref{2lin}), the push-forward of $\pi$ by this map has necessarily
a cross-term. It is therefore necessary to admit an additional braiding
(nontrivial cross-relations) between links. The full Poisson
structure on two links will be then given by
\be
\pi _{++}(a,b)= \pi _+(a) \oplus \pi _+(b)\oplus\pib
(a,b),\qquad a,b\in G, 
\ee
according to the obvious decomposition of the tangent space
$T_{(a,b)}(G\times G)$. The cross-term $\pib$ can be always written
in the form
\be
\pib (a,b)= (a,e)\phi (a,b)(e,b),
\ee
where
\be\label{is}  
\phi (a,b) \in \gotG _{(1)}\bigwedge \gotG _{(2)}\cong
\gotG \otimes \gotG 
\ee
($e\in G$ is the group unit).

Let us mention that the braiding in the corresponding quantum case  
will be given by a cross-product coalgebra structure on
$H\otimes H$:
\be
\Dr _{++} = (\id\otimes Q\otimes \id )(\Dr _{+}\otimes \Dr _{+}),
\ee
where $Q\colon H\otimes H\to H\otimes H$ is a linear invertible
map. 

Now we state our results (their proof is given in the Appendix).
\begin{Prop}\label{tw1}
The map given by formula (\ref{2lin}) is Poisson if and only if
\be\label{pier}
\phi (a,b) = -r +\psi (a,b),\qquad \mbox{\rm where}\;\;\;\;
\psi (ag^{-1},gb)=\Ad _g \psi (a,b) .
\ee
\end{Prop}
\begin{Rem}
{\rm The above condition for $\psi$ has a general solution given
by} \ $\psi (a,b)=\Ad _b f(ab)$, {\rm where $f$ is an arbitrary 
function on $G$ with values in $\gotG \otimes \gotG$.}
\end{Rem}
\begin{Prop} ({\bf constant case}). \ Assume that $\phi = const
= -w$ satisfies (\ref{pier}). Then $w-r$ is invariant and we may
assume that it is symmetric (one can modify $r$, if necessary). 
The bivector field $\pi _{++}$ is Poisson if and only if
\be\label{cybe}
[w,w]  \equiv [w_{12},w_{13}]+[w_{12},w_{23}]+[w_{13},w_{23}]=0.
\ee
\end{Prop}
The above Proposition clearly states that in the constant case, 
our condition is equivalent to the quasitriangularity.

\section{Composability}

Another natural requirement concerning two links is that the composition 
(of parallel transports) should be a Poisson map, more precisely,
the multiplication
\be
G\times G \ni (a,b)\mapsto ab\in G
\ee
should be a Poisson map as a map from $(G\times G, \pi _{++} )$ 
to $(G,\pi _+ )$. 
\begin{Prop} The multiplication is a Poisson map if and only if 
\be
\phi (a,b)_{\rm anti} = -r
\ee
(`{\rm anti}' refers to the antisymmetric part).
\end{Prop}
\dowod
For $\phi (a,b)=X\otimes Y\cong X_{(1)}\w Y_{(2)}$ we have $\pib
(a,b)= (aX)_{(1)}\w (Yb)_{(2)}$ which is mapped by the
multiplication on 
$$(aXb)\w (aYb)=a(X\w Y)b=a (\phi (a,b)-P\phi (a,b))b=a(2\phi
(a,b)_{\rm anti})b .$$
The result being valid for arbitrary $\phi$ (by linearity), the
map $(a,b)\mapsto ab$ is Poisson if and only if
$$\pi _{++}(a,b)=(ra+ar)b+a(rb+br)+2a(\phi (a,b)_{\rm anti})b =
rab + abr,$$
which yields the required result.

\qed

The above result complements Proposition~\ref{tw1}: the
antisymmetric part must be constant.

Concerning the quantum case, we have the following proposition.
\begin{Prop}
The multiplication map $m\colon H\otimes H\to H$ is a coalgebra
homomorphism as a map from $(H\otimes H,\Dr _{++})$ to
$(H,\Dr_{+})$,
\be\label{homo} 
\Dr _{+}m= (m\otimes m)\Dr _{++},
\ee
if and only if $Q$ is given by
\be\label{Q}
Q=P(m\otimes m)(\id\otimes P\sfR _D),
\ee
where $\sfR _D := \sfR ^{-1}$ (which is the $R$-matrix of
Drinfeld's type). 
$\Dr _{++}$ is coassociative if and only if $\sfR _D$  satisfies
\be\label{quasi}
(\Dr \otimes \id )\sfR _D = (\sfR _D)_{13}(\sfR _D)_{23},\qquad
(\id \otimes \Dr )\sfR _D = (\sfR _D)_{12}(\sfR _D)_{13},
\ee
i.e. the quasitriangularity property.
\end{Prop}
\dowod
First we show that (\ref{Q}) satisfies (\ref{homo}). Using $\Dr
_{+}(a)= \Dr (a)\, \sfR $, we have
$$ (m\otimes \id)(\id \otimes Q)(\Dr _{+}\otimes \id )=
(m\otimes \id)(\Dr \otimes \id )
$$
(diagrams!), hence 
$$ (m\otimes m)(\id\otimes Q\otimes \id )(\Dr _{+}\otimes \Dr
_{+})= (m\otimes m)(\id \otimes P\otimes \id )(\Dr \otimes \Dr _{+})=
\Dr _{+}m.
$$
We shall show that this solution is unique. It is sufficient to
show that 
\be\label{zero}
m(X\otimes \id )\Dr _{+}=0
\ee
implies $X=0$ for each $X\in \End (H)$. 
 Assuming (\ref{zero}), and using
\be
(\id\otimes m)(\Dr _{+}\otimes S)P\Dr = (m\otimes \id )
(\id\otimes P)(\sfR \otimes \id )
\ee
($S$ is the antipode), we have 
$$
0=m(m\otimes \id )(X\otimes \id \otimes S)(\Dr _{+}\otimes \id
)P\Dr = m (X\otimes \id )(\id\otimes m)(\Dr _{+}\otimes S)P\Dr 
$$
\be
=m(X\otimes \id )(m\otimes \id )(\id \otimes P)(\sfR \otimes \id ).
\ee
Denoting the last expression by $F$, we have
$$
0=m(F\otimes \id )(m\otimes \id )(\id \otimes P)(\sfR _D \otimes
\id )=X. 
$$
For the second part of the proof, recall that $\Dr _{++}$ is
coassociative if and only if 
\be\label{dwawar}
(\Dr _+\otimes \id )Q=(\id \otimes Q)(Q\otimes \id )(\id
\otimes\Dr _+),\;\;\;\;
(\id\otimes\Dr _+ )Q=(Q\otimes\id )(\id\otimes Q)(\Dr
_+\otimes\id ).
\ee
We denote by $P_{k_1k_2\ldots k_n}$ the permutation of the
tensor product sending the $j$-th tensor factor to the $k_j$-th place.
Let $L$ and $R$ denote the left and right hand side of the first 
condition in (\ref{quasi}). This condition is equivalent to
$A=B$, where
\be
A = (m\otimes m\otimes\id )P_{13524}(L\otimes\id\otimes\id ),\qquad
B = (m\otimes m\otimes\id )P_{13524}(R\otimes\id\otimes\id ).
\ee
Now note that the first condition in (\ref{dwawar}) can be
written as
\be
(\id\otimes \id\otimes m)P_{3124}(\id\otimes A)(\id\otimes\Dr _+)=
(\id\otimes \id\otimes m)P_{3124}(\id\otimes B)(\id\otimes\Dr _+).
\ee
or,
\be\label{war1}
A\Dr _+= B\Dr _+.
\ee
Therefore, the first condition in (\ref{quasi}) implies the
first condition in (\ref{dwawar}). The converse implication
follows by applying (\ref{war1}) to $I$, and using $L=A(I\otimes
I)$, $R=B(I\otimes I)$.

\qed

Using (\ref{Q}), we can derive a formula for the braiding
on the level of matrix quantum groups. Calculating the transposed $Q$,
\be
Q^T=(\id \otimes \sfR _D^TP\otimes \id )(\Dr\otimes\Dr )P
\ee
($\Dr=m^T$ is the comultiplication in the `function algebra':
$\Dr v^a{_l}= v^a{_j}\otimes v^j{_l}$)
on $w^a{_l}\otimes v^k{_b}$, we obtain
\be\label{QT}
Q^T(w^a{_l}\otimes v^k{_b})=(\id\otimes \sfR _D^T\otimes \id )
(v^k{_j}\otimes w^a{_m}\otimes v^j{_b}\otimes w^m{_l})
= R^{aj}{_{mb}}v^k{_j}\otimes w^m{_l},
\ee
where $R^{aj}{_{mb}} = \left\langle R_D^T ,w^a{_m}\otimes
v^j{_b}\right\rangle $ is the $R$-matrix  of FRT type
($R=P\hat{R}$), 
hence the crossed multiplication rule in the tensor product of
two copies of algebras with the `plus' structure
\be 
\hat{R}v\Ttimes v = v\Ttimes v \hat{R}_{21},\qquad 
\hat{R}w\Ttimes w = w\Ttimes w \hat{R}_{21},
\ee
is given by
\be\label{braid} 
 w^a{_l}v^k{_b}= v^k{_s}R^{aj}{_{mb}}w^t{_l}.
\ee
One can see directly that the composition is a homomorphism, by
checking that $v^j{_a}w^a{_l}$ satisfy again the `plus' relations:
$$
\hat{R}^{pq}{_{jk}}v^j{_a}w^a{_l}v^k{_b} w^b{_m}=
\hat{R}^{pq}{_{jk}}v^j{_a}v^k{_s}\hat{R}^{sa}{_{tb}}w^t{_l} w^b{_m}=
v^p{_j}v^q{_k}\hat{R} ^{kj}{_{sa}}
w^s{_t} w^a{_b}\hat{R} ^{bt}{_{ml}}
$$
$$
= v^p{_j}w^j{_t}v^q{_a}w^a{_b}\hat{R} ^{bt}{_{ml}}.
$$
(One can also see that (\ref{QT}) and the transpose of (\ref{dwawar})
define $Q^T$ consistently on higher order polynomials, due to
Yang-Baxter equation.)

\section{Final remarks}

We have seen that in order to implement gauge transformations
(\ref{1lin}), (\ref{2lin}) as Poisson (or quantum)
homomorphisms, the Poisson (quantum) gauge group should be
quasitriangular, in the {\em real sense}. This excludes in
particular the standard deformations of compact simple groups,
which are always quasitriangular in the {\em imaginary sense}
(cf.~also \cite{stand}). 

A simple imaginary quasitriangular case is provided by the standard
deformation of $SU(2)$. Expanding the standard $R$-matrix for $SU(2)$ up
to linear terms in $\er =q-1$, we get
$$
R= \left( \ba{cccc} q^{\frac12}  & 0 & 0 & 0 \\
                       0 & q^{-\frac12} & 0 & 0 \\
  0 &  q^{\frac12}- q^{-\frac32} & q^{-\frac12} & 0 \\
                       0 & 0 & 0 & q^{\frac12} \ea\right)
\sim \; I + \er 
\left( \ba{cccc} \frac12  & 0 & 0 & 0 \\
                       0 & -\frac12 & 0 & 0 \\
                       0 & 2 & -\frac12 & 0 \\
                       0 & 0 & 0 & \frac12 \ea\right)
$$
$$
= I + \er (2X_-\otimes X_+ +\frac12 \sir _3\otimes \sir _3 )
= I + iw,$$
where
$$ X_+ := \left(\ba{cc} 0 & 1 \\ 0 & 0 \ea\right) ,\qquad
X_- := \left(\ba{cc} 0 & 0 \\ 1 & 0 \ea\right) ,
$$
and
$$iw = \er X_-\wedge X_+ + \er (X_-\otimes X_+ +X_+\otimes X_-+
\frac12 \sir _3\otimes \sir _3 )
$$
$$=
\frac{i\er}{2}\sir _1\wedge\sir _2 + \frac{\er}{2}
( \sir _1\otimes \sir _1 + \sir _2\otimes \sir _2 + 
\sir _3\otimes \sir _3)= i (r - i s),$$
$$ r := \frac{\er}{2}\sir _1\wedge\sir _2 \in \Ldwa su(2),\qquad
s:= \frac12 \sum_{j=1}^3 \sir _j\otimes\sir _j \in 
\Odwa _{\rm symm} su(2).$$ 
(Here $\sigma _j$ are the Pauli matrices.) 
We see that the symmetric part of $w$ is imaginary.

Consider now the case of the (standard) quantum
$SU(N)$. We assume therefore that the matrices $v$ and $w$ in
(\ref{braid}) are unitary. It is easy to see that the conjugation
of (\ref{braid}) gives again (\ref{braid}) (does not lead to new
relations) if and only if $\hat{R}$ is unitary. This is in
contradiction with the fact that it is actually self-adjoint
(but not involutive).

Finally, we remark that having transformations (\ref{1lin}),
(\ref{2lin}) implemented on the Poisson (quantum) level, they are
automatically implemented on the whole one-dimensional lattice.

\section{Appendix}

\subsection{Proof of Prop.~1.1}

Let $\Phi $ be the map in (\ref{2lin}) and let us denote
isomorphism (\ref{is}) explicitly by
$$ \gotG\otimes \gotG \ni X\otimes Y \mapsto 
(X\otimes Y)_{[12]}:= X_{(1)}\w
Y_{(2)}\in \gotG _{(1)}\bigwedge \gotG _{(2)}.
$$
For $X,Y\in \gotG$, we have $ \Phi _* (Xg) = 
(-ag^{-1}X)_{(1)} + (Xgb)_{(2)}$ and
$$ \Phi _* (Xg)\w \Phi _* (Yg)=(ag^{-1}X)_{(1)}\w (ag^{-1}Y)_{(1)} 
+ (Xgb)_{(2)}\w (Ygb)_{(2)} - 
$$
$$
-(ag^{-1}X)_{(1)}\w (Ygb)_{(2)} -(Xgb)_{(2)}\w (ag^{-1}Y)_{(1)} ,
$$
hence, for $r\in \Ldwa \gotG $,
\be\label{Prg}
\Phi _* (rg) = (ag^{-1}r)_{(1)} + (rgb)_{(2)} -
(ag^{-1})_{(1)}r_{[12]}(gb)_{(2)} .
\ee
Similarly, using
$$ \Phi _* (gX)=(-aXg^{-1})_{(1)}+(gXb)_{(2)} $$
we obtain
\be
\Phi _* (gr) = (arg^{-1})_{(1)} + (grb)_{(2)} -
(a,g)r_{[12]}(g^{-1},b).
\ee
Since $a_{(1)} (X\otimes Y)_{[12]} b_{(2)}=(aX)_{(1)}\w
(Yb)_{(2)}$, and
\be
\Phi _* (aX)_{(1)}\w \Phi _* (Yb)_{(2)}= (aXg^{-1})_{(1)}\w
(gYb)_{(2)} = (a,g)X_{(1)}\w Y_{(2)} (g^{-1},b),
\ee
we have
\be
\Phi _* \pib (a,b) =  (a,g)\phi (a,b)_{[12]} (g^{-1},b).
\ee
Of course
\be\label{Ppi+}
\Phi _* \pi _+ (a)= ((ra+ar)g^{-1})_{(1)},\qquad 
\Phi _* \pi _+ (b)= (g(ra+ar))_{(2)}.
\ee
Inserting (\ref{Prg})--(\ref{Ppi+}) into
$$ \Phi _* (\pi _+(a) +\pi _+(b) +\pib (a,b) +\pi (g))=
(rag{\!}^{-1} +ag{\!}^{-1}r)_{(1)} + 
(rgb + gbr)_{(2)} +\pib (ag^{-1},gb)
$$
we obtain 
\be
\phi (ag^{-1},gb)+r=\Ad _g (\phi (a,b) +r).
\ee

\subsection{Proof of Prop.~1.3}

For $X\in \gotG$ we set $X^L(g):=gX$, $X^R(g):=Xg$ and,
similarly for higher rank tensors. In particular,
$\pi _+= r^L + r ^R$.
Expanding $w$ in some basis $X_k$ in $\gotG$, we have
$$-\pib (a,b)=a_{(1)}w_{[12]}b_{(2)}= 
a_{(1)}w^{kl}X_{k(1)}\w X_{l(2)}b_{(2)}=
w^{kl}X_{k(1)}^L\w X_{l(2)}^R=: w_{[12]}^{LR}.
$$
Since $\pi _{++}=
r^L_{(1)}+r^R_{(1)}+r^L_{(2)}+r^R_{(2)}-w_{[12]}^{LR}$, we have
$$
[\pi _{++},\pi _{++}]= -2
[r^L_{(1)}+r^R_{(1)}+r^L_{(2)}+r^R_{(2)},w_{[12]}^{LR}]
+ [ w_{[12]}^{LR},w_{[12]}^{LR}]
$$
$$
=-2[r^L_{(1)}+r^R_{(2)},w_{[12]}^{LR}]+ [ w_{[12]}^{LR},w_{[12]}^{LR}].
$$
We set $r=\frac12 r^{ij}X_i\w X_j$ and calculate the first term:
$$
[r^L_{(1)},w^{LR}_{[12]}]=
\frac12 r^{ij}w^{kl}[X^L_{i(1)}\w X^L_{j(1)},X^L_{k(1)}\w X^R_{l(2)}]
=r^{ij}w^{kl}X^L_{i(1)}\w [X_j,X_k]^L_{(1)}\w X^R_{l(2)}
$$
$$
\cong 
r^{ij}w^{kl}X^L_{i(1)}\w [X_j,X_k]^L_{(1)}\otimes X^R_{l(2)}
$$
$$
=r^{ij}w^{kl}(X^L_{i(1)}\otimes [X_j,X_k]^L_{(1)}-
 {}[X_j,X_k]^L_{(1)}\otimes X^L_{i(1)})\otimes X^R_{l(2)}
= [r_{12},w_{23}]+[r_{12},w_{13}],
$$
the last expression being interpreted as an element of $\gotG
_{[11]2}:= (\Ldwa \gotG )_{(1)}^L\w \gotG _{(2)}^R$.
Similarly,
\be
 {}[r^R_{(2)},w_{[12]}^{LR}]= [r_{23},w_{13}]+[r_{23},w_{12}],
\ee
where the right hand side is understood as an element of 
$\gotG _{1[22]}:=\gotG ^L_{(1)} \w (\Ldwa \gotG )^R _{(2)}$,
and also
\be
 {}[w_{[12]}^{LR},w_{[12]}^{LR}]= 2[w_{12},w_{13}]
- 2[w_{13},w_{23}],
\ee
where the terms on the right hand side are elements of
$\gotG _{1[22]}$ and $\gotG _{[11]2}$, respectively.
It follows that $\pi _{++}$ is Poisson if and only if
\be\label{cybe1}
 {}[r_{12},w_{13}]+ [r_{12},w_{23}]+[w_{13},w_{23}] =0
\ee
and
\be\label{cybe2}
 {}[w_{12},w_{13}]+ [w_{12},r_{23}]+[w_{13},r_{23}] =0.
\ee
Due to the invariance of the symmetric part $w-r$ of $w$, both
equations are equivalent to (\ref{cybe}).

\end{document}